\documentclass[12pt]{article} 
\usepackage[sectionbib]{natbib}
\usepackage{array,epsfig,fancyheadings,rotating}
\usepackage[]{hyperref}  
\usepackage{sectsty, secdot}
\usepackage{xcolor}
\sectionfont{\fontsize{12}{14pt plus.8pt minus .6pt}\selectfont}
\renewcommand{\theequation}{\thesection\arabic{equation}}
\subsectionfont{\fontsize{12}{14pt plus.8pt minus .6pt}\selectfont}

\usepackage{comment}
\usepackage{amsmath}
\usepackage{amssymb}
\usepackage{amsfonts}
\usepackage{multirow}
\usepackage{amsthm}
\usepackage[mathscr]{eucal}
\usepackage{cite}
\usepackage{algorithm}
\usepackage{algorithmic}
\usepackage{setspace}
\usepackage{svg}
\usepackage{subfigure}
\usepackage{graphicx}

\makeatletter
\renewcommand{\ALG@name}{Algorithm}
\newenvironment{breakablealgorithm}
  {
   \begin{flushleft}
     \refstepcounter{algorithm}
     \hrule height.8pt depth0pt \kern2pt
     \renewcommand{\caption}[2][\relax]{
       {\raggedright\textbf{\ALG@name~\thealgorithm} ##2\par}%
       \ifx\relax##1\relax 
         \addcontentsline{loa}{algorithm}{\protect\numberline{\thealgorithm}##2}%
       \else 
         \addcontentsline{loa}{algorithm}{\protect\numberline{\thealgorithm}##1}%
       \fi
       \kern2pt\hrule\kern2pt
     }
  }{
     \kern2pt\hrule\relax
   \end{flushleft}
  }
\makeatother

\newtheorem{corollary}{Corollary}

\theoremstyle{plain}

\newtheorem{assumption}{Assumption}
\newtheorem{thm}{Theorem}
\begin{document}

\renewcommand{\baselinestretch}{1}

\markright{ \hbox{\footnotesize\rm 
}\hfill\\[-13pt]
\hbox{\footnotesize\rm
}\hfill }

\markboth{\hfill{\footnotesize\rm FIRSTNAME1 LASTNAME1 AND FIRSTNAME2 LASTNAME2} \hfill}
{\hfill {\footnotesize\rm FILL IN A SHORT RUNNING TITLE} \hfill}

\renewcommand{\thefootnote}{}
$\ $\par


\fontsize{12}{14pt plus.8pt minus .6pt}\selectfont \vspace{0.8pc}
\centerline{\large\bf Optimal Subsampling for}
\vspace{2pt} 
\centerline{\large\bf Large Scale Elastic-Net Regression}
\vspace{.4cm} 
\centerline{Hang Yu$^{1}$, Zhenxing Dou$^{1}$, Zhiwei Chen$^{1}$ and Xiaomeng Yan$^{2}$}
\vspace{.4cm} 
\centerline{\it $^{1}$Beihang University and $^{2}$Capital University of Economics and Business}
\vspace{.55cm} \fontsize{9}{11.5pt plus.8pt minus.6pt}\selectfont
\begin{quotation}
\noindent {\it Abstract:}
Datasets with sheer volume have been generated from fields including computer vision, medical imageology, and astronomy whose large-scale and high-dimensional properties hamper the implementation of classical statistical models. To tackle the computational challenges, one of the efficient approaches is subsampling which draws subsamples from the original large datasets according to a carefully-design task-specific probability distribution to form an informative sketch. The computation cost is reduced by applying the original algorithm to the substantially smaller sketch. Previous studies associated with subsampling focused on non-regularized regression from the computational efficiency and theoretical guarantee perspectives, such as ordinary least square regression and logistic regression. In this article, we introduce a randomized algorithm under the subsampling scheme for the Elastic-net regression which gives novel insights into $L_1$-norm regularized regression problem. To effectively conduct consistency analysis, a smooth approximation technique based on $\alpha$-absolute function is firstly employed and theoretically verified. The concentration bounds and asymptotic normality for the proposed randomized algorithm are then established under mild conditions. Moreover, an optimal subsampling probability is constructed according to A-optimality. The effectiveness of the proposed algorithm is demonstrated upon synthetic and real data datasets. 
\par

\vspace{9pt}
\par
\end{quotation}\par


\def\thefigure{\arabic{figure}}
\def\thetable{\arabic{table}}

\renewcommand{\theequation}{\thesection.\arabic{equation}}

\fontsize{12}{14pt plus.8pt minus .6pt}\selectfont

\section{ Introduction}
\par
The evolution of information and technology brings a rapid growth of big data in scientific research fields ranging from social networking to medical imaging. Dealing with massive data whose sample size and dimensional features increase out of the capacity is always a daunting challenge for classical statistical models. There have been extensive studies on randomization-based methods to cut down the data quantity utilized in the model computing process and thereby reduce the high computation costs, among which, benefited from its easy-to-use property, subsampling has gained popularity, including matrix computation and decomposition~\citep[e.g.,][]{drineas2006fast1,drineas2006fast2,drineas2006fast3}, Markov Chain Monte Carlo algorithms~\citep[e.g.,][]{lin2000noisy,beaumont2003estimation,andrieu2009pseudo,chazal2015subsampling,dang2019hamiltonian}, and deep learning~\citep[e.g.,][]{huijben2020deep,van2021active}.

The mechanism of subsampling is to construct a substantially smaller representative sketch by drawing subsamples from the given large-scale data at hand in accordance with some carefully-designed subsampling probability (SSP) and approach the machine learning problem based on the sketch. To retain the most relevant information in the sketch and, furthermore, keep the result computed from the sketch close to the result from the original dataset, the informative samples are supposed to be assigned with larger probabilities. Since the seminal works of \citet{drineas2011faster, drineas2012fast} that facilitated leverage matrix to construct basic leveraging (BLEV) and approximate leveraging subsampling probability to
speed up the least squares approximation, the idea of subsampling has been successfully applied to regression models.
\citet{zhu2015optimal} proposed optimal subsampling and predictor-length subsampling for ordinary least squares regression and developed theoretical guarantees for the subsampled estimator in terms of the consistency to the full sample estimator. For large-scale logistic regression, \citet{fithian2014local} introduced local-case control sampling to deal with imbalanced data sets, and an optimal subsampling probability was designed by \citet{wang2018optimal} to minimize the trace of the asymptotic variance-covariance matrix of the subsampled estimator, i.e., using A-optimality \citep{kiefer1959optimum} as an optimal criterion. By extending the work of \citet{wang2018optimal}, \citet{yao2019optimal} designed optimal subsampling probability with A-optimality for softmax regression, which is also called multinomial logistic regression. L-optimality is adopted as a reference for designing optimal SSP for quantile regression \citep{wang2021optimal}. Another popular randomization-based method to make statistical models scalable to large datasets is random projection, see \citet{rokhlin2008fast}, \citet{dhillon2013new}, \citet{iyer2016randomized}, \citet{drineas2016randnla} and \citet{clarkson2017low}. 

\par

However, besides the quantity of data, as the dimensionality of the data escalates towards higher orders of magnitude and collinearity appears as a major concern for statistical inference, classical ordinary least squares regression and ridge regression are inadequate due to their limited interpretability and robustness.
Regression models equipped with {\it $L_1$}-norm penalty term such as Lasso \citep{tibshirani1996regression} and Elastic-net \citep{zou2005regularization} have been demonstrated with powerful feature selection capabilities to address high-dimensional data problems by forcing some coefficients to be exactly zeros. To further elaborate, Elastic-net regression combines {\it $L_1$}-norm and {\it $L_2$}-norm penalty terms to balance model stability and feature selection capability, which contribute to its wide utilization in down-to-earth engineering problems \citep[e.g.,][]{zou2003regression,li2022estimating}. \citet{zhou2015reduction} reveals the relationship between Elastic-net and support vector machine and designs a parallel solver for the Elastic-net to accelerate GPU computations. The process of attacking deep neural networks (DNNs) via adversarial examples was formulated as an Elastic-net regularized optimization problem, see \citet{chen2018ead}. The Elastic-net regularized Cox proportional hazards model was proposed in machine learning algorithm to predict the overall survival among re-infected COVID-19 patients\citep{ebrahimi2022predictive}.
\par
A-optimality comes from the statistical motivation that minimizing the asymptotic MSE of the parameters. In subsampling algorithm, minimizing the asymptotic MSE implies that as the number of subsamples increases, the subsample estimator converges to the full sample estimator at a faster rate. It is noteworthy that the aforementioned A-optimality-based subsampling methods share a common property that the criterion functions are infinitely differentiable, leading to relatively straightforward consistency analysis and asymptotic normality derivation by utilizing Taylor series expansion. 
Nonetheless, the non-differentiability of {\it $L_1$}-norm penalty term involved in the Elastic-net criterion function precludes the routine asymptotic analysis on the Elastic-net regression estimator, 
as a result of which, designing optimal subsampling probability (OSP) based on A-optimality for Elastic-net regression to reduce the burden on computing resources and verify its validation based on asymptotic analysis is necessary but difficult. 

In this paper, we advocate an effective randomized algorithm based on subsampling technique for Elastic-net regression, and our main methodological and theoretical contributions are highlighted as follows,
\begin{itemize}
 \setlength\itemsep{-0.25em}
    \item[1.] {\it Smooth approximation to Elastic-net criterion function.} The {\it $L_1$}-norm regularization term introduces a non-differentiable property to the criterion function of Elastic-net, thereby limiting the feasibility of asymptotic analysis on the estimator. We present an innovative smooth approximation for Elastic-net regression by utilizing the $\alpha$-absolute function. Theoretically, we prove that the global unique minimum of smooth Elastic-net converges to the original one as $\alpha\to\infty$. Such smooth approximation provides convenience for subsequent theoretical analysis and numerical optimization algorithm design.
    \item[2.] {\it Optimal subsampling design.} In this paper, we rigorously prove the consistency and asymptotic normality of the estimator from a general subasampling algorithm for smooth Elastic-net regression. Additionally, according to the asymptotic distribution, we derived the optimal subsampling probability which minimizes the asymptotic MSE of the subsampling estimator. Since such subsampling probability may not exist, we have designed a pseudo-optimal subsampling probability as a substitution.
    \item[3.] {\it A two-step algorithm.} Considering that the pseudo-optimal subsampling probability depends on the full sample set, we propose a two-step algorithm to approximate the optimal subsampling procedure. In the first step, the pseudo-optimal subsampling probability is calculated on a small subsample set via uniform subsampling. The second step is to sample from the full sample set by using the pseudo-optimal subsampling probability from the first step. We prove consistency and asymptotic normality of the estimator from the two-step algorithm.
\end{itemize}
\par
The rest this article is organized as follows. In \S\ref{sec:smooth}, We provide an overview of the Elastic-net regression and establish a smooth approximation to the original criterion function. The convergence of the global minimum of the smooth Elastic-net to the original problem is developed. In \S\ref{sec:mainres}, we establish consistency and asymptotic normality for the resultant estimator of the general subsampling procedure. In \S\ref{sec:optimal}, based on the asymptotic MSE, we design pseudo-optimal SSP for smoothed Elastic-net regression in accordance with A-optimality. In \S\ref{sec:twostep}, we develop a two-step algorithm to implement the randomized Elastic-net regression efficiently. Moreover, consistency and asymptotic normality of the resultant estimator obtained from the two-step algorithm are also established. In \S\ref{sec:simulation}, we design four examples to examine the theoretical results about Algorithm \ref{alg:2step-alg} proposed in \S\ref{sec:twostep}. In \S\ref{sec:realdata}, the performance of the two-step algorithm is elaborated through two practical applications.

\section{Problem Setup}\label{sec:smooth}

Given samples $\left\{(\boldsymbol{x}_1,y_1),(\boldsymbol{x}_2,y_2),\dots,(\boldsymbol{x}_N,y_N)\right\}$ with covariates $\boldsymbol{x}_n\in \mathbb{R}^p$ and response $y_n\in \mathbb{R}$, they are identically independently distributed as random variable pair $(\boldsymbol{x}, y)\in\mathbb{R}^{p+1}$ and the conditional expectation of $y$ given $\boldsymbol{x}$ satisfies the linear equation $\mu = \mathbb{E}(y|\boldsymbol{x}) = \boldsymbol{x}^\top\boldsymbol{\beta}$. Let $\mathbf{Y} = \left (y_1,\dots,y_N\right)^\top$ and $\mathbf{X} = \left(\boldsymbol{x}_1,\dots,\boldsymbol{x}_N\right)^\top$, the relationship between $\mathbf{X}$ and $\mathbf{Y}$ is established through a linear model as below,
\begin{equation}\label{GaussianLinearModel}
\mathbf{Y} = \mathbf{X}\boldsymbol{\beta} + \sigma \boldsymbol{\epsilon},
\end{equation}
where $\boldsymbol{\beta}$ is the coefficient vector and $\boldsymbol{\epsilon} = (\epsilon_1,...,\epsilon_N)^{\top}$ is white noise with zero expectation and unit variance. All subsequent theoretical analyses in this paper will be conducted under the aforementioned assumptions, and we will not elaborate on it further. To solve for the coefficient $\boldsymbol{\beta}$ with desired properties, the Elastic-net estimate with regularized hyper-parameter $\lambda$ is obtained by minimizing the following criterion function, 
\begin{equation}\label{ElasticNet}
\mathcal{L}(\boldsymbol{\beta}; \lambda ,\eta) = \frac{1}{2}
\|\mathbf{Y} - \mathbf{X}\boldsymbol{\beta}\|^2 + \lambda \mathcal{P}_{\eta}(\boldsymbol{\beta}), 
\end{equation}
where the norm $\| \cdot \|$ is {\it $L_2$}-norm, and the penalty function with $\eta\in (0,1)$ is defined by
\begin{equation*}
\mathcal{P}_{\eta}\left(\boldsymbol{\beta}\right) =  \frac{1}{2}(1-\eta){\|\boldsymbol{\beta}\|}^2_2 +  \eta{\|\boldsymbol{\beta}\|}_1.
\end{equation*}
In practice, $\lambda > 0$ and $\eta \in (0,1)$ are commonly chosen by cross-validation with the purpose of model scalability and predictive performance. 
We represent the unknown optimal estimator of problem \eqref{ElasticNet} as $\widehat{\boldsymbol{\beta}}_{opt}$ such that
\begin{equation*}
\widehat{\boldsymbol{\beta}}_{opt} = \mathrm{arg}\min\limits_{\boldsymbol{\beta}}\mathcal{L}(\boldsymbol{\beta};\lambda ,\eta).
\end{equation*}
\par
\textbf{Smooth Approximation}. In order to facilitate the subsequent numerical optimization and asymptotic normality analysis, $\alpha$-absolute function \citep{schmidt2007fast} is introduced to approximate $L_1$ term smoothly. The smooth approximation to the absolute function $|x|$ in terms of the $\alpha$-absolute function is defined as
\begin{equation*}
{|x|}_{\alpha} := \frac{1}{\alpha} \left[ \log\left(1+\exp\left(-\alpha x\right)\right) + \log\left(1+\exp\left(\alpha x\right)\right) \right],
\end{equation*}
which is equipped with the following properties:
\begin{equation*}
\nabla {|x|}_{\alpha} = \left(1 + \exp\left(-\alpha x\right)\right)^{-1} - \left(1 + \exp\left(\alpha x\right)\right)^{-1},
\end{equation*}
\begin{equation}\label{alphaAbsH}
{\nabla}^{2}{|x|}_{\alpha} = 2\alpha \exp\left(\alpha x\right) / \left(1+\exp\left(\alpha x\right)\right)^2,
\end{equation}
\begin{equation*}
\lim_{\alpha \rightarrow + \infty} {|x|}_{\alpha} = |x|.
\end{equation*}

Therefore, $\mathcal{L}_{\alpha}(\boldsymbol{\beta};\lambda ,\eta)$, which is constructed based on $\alpha$-absolute function, is a smooth approximation to $\mathcal{L}(\boldsymbol{\beta};\lambda ,\eta)$.
\begin{equation}
\label{smoothElasticNet}
\mathcal{L}_{\alpha}(\boldsymbol{\beta};\lambda ,\eta)  = \frac{1}{2}\|\mathbf{Y} - \mathbf{X}\boldsymbol{\beta}\|^2 + \lambda \left(\frac{1}{2}(1-\eta){\|\boldsymbol{\beta}\|}^2_2 + \eta{\|\boldsymbol{\beta}\|}_{\alpha}\right),
\end{equation}
where
$${\|\boldsymbol{\beta}\|}_{\alpha}  := \sum_{j=1}^p |\boldsymbol{\beta}_j|_{\alpha}. $$
Note that ${\|\boldsymbol{\beta}\|}_{\alpha}$ operates smooth approximation towards vector $\boldsymbol{\beta}$ elementwisely, which is well-defined to reflect the separable nature of $L_1$ norm. Without otherwise specification, the following numerical and theoretical analysis towards the Elastic-net estimator is based on the smoothed criterion function, i.e., $\mathcal{L}_{\alpha}(\boldsymbol{\beta};\lambda ,\eta)$. Correspondingly, the unknown optimal estimator of the smoothly approximated Elastic-net criterion can be computed by
\begin{equation}
\label{betaOsa}
\widehat{\boldsymbol{\beta}}_{osa} = \mathrm{arg}\min \limits_{\boldsymbol{\beta}}\mathcal{L}_{\alpha}(\boldsymbol{\beta};\lambda ,\eta).
\end{equation}
\par
In practical implementation of the algorithm, $\alpha$ is often chosen from the positive integer sequence $\{1, 2, 3, ..., m, ...\}$.
\begin{thm}\label{thm1}
Let $m$ be a positive integer, $\mathcal{L}_{m}(\boldsymbol{\beta};\lambda ,\eta)$ be the corresponding smooth approximation to Elastic-net criterion function $\mathcal{L}(\boldsymbol{\beta};\lambda ,\eta)$ and $\widehat{\boldsymbol{\beta}}_{osa}^{(m)}$ be the unknown optimal estimator of $\mathcal{L}_{m}(\boldsymbol{\beta};\lambda ,\eta)$. If the covariate dimension $p$ is fixed, $\widehat{\boldsymbol{\beta}}_{osa}^{(m)}$ converges to $\widehat{\boldsymbol{\beta}}_{opt}$ as $m\to\infty$, i.e.,
\begin{equation}
\lim_{m \to \infty}\widehat{\boldsymbol{\beta}}_{osa}^{(m)} = \widehat{\boldsymbol{\beta}}_{opt}.
\end{equation}
\end{thm}
Theorem \ref{thm1} indicates that if the hyperparameter $\alpha$ in \eqref{smoothElasticNet} is chosen from the positive integer sequence $\{1, 2, 3, ..., m, ...\}$, then the unknown optimal estimator of \eqref{smoothElasticNet} converges to the unknown optimal estimator of \eqref{ElasticNet} as $\alpha\to\infty$.
\par
\textbf{Newton's Method}. $\mathcal{L}_{\alpha}(\boldsymbol{\beta};\lambda ,\eta)$ is infinitely differentiable and strictly convex when $\eta \in (0,1)$. Analytically, there is no explicit closed-form solution to the $\widehat{\boldsymbol{\beta}}_{osa}$ and a commonly used optimization tool is Newton's method. Specifically, for smooth Elastic-net regression, Newton's method iteratively applies the following formula until $\widehat{\boldsymbol{\beta}}_{osa}^{(t)}$ converges in terms of given tolerance,
\begin{equation}\label{NewtonMethod}
\widehat{\boldsymbol{\beta}}_{osa}^{(t+1)} = \widehat{\boldsymbol{\beta}}_{osa}^{(t)} - \biggl[\frac{\partial^2 \mathcal{L}_{\alpha}(\boldsymbol{\beta};\lambda ,\eta)}{\partial \boldsymbol{\beta} ^2}\big|_{\boldsymbol{\beta} = \widehat{\boldsymbol{\beta}}_{osa}^{(t)}}\biggl]^{-1}\frac{\partial \mathcal{L}_{\alpha}(\boldsymbol{\beta};\lambda ,\eta)}{\partial \boldsymbol{\beta}}\big|_{\boldsymbol{\beta} = \widehat{\boldsymbol{\beta}}_{osa}^{(t)}},
\end{equation}
in which
\begin{equation}\label{L_alpha_G}
\frac{\partial \mathcal{L}_{\alpha}(\boldsymbol{\beta};\lambda ,\eta)}{\partial \boldsymbol{\beta}} = \sum_{n=1}^N(\boldsymbol{\beta}^{T} \boldsymbol{x}_n - y_n)\boldsymbol{x}_n + \lambda(1-\eta)\boldsymbol{\beta} + \lambda\eta \boldsymbol{g}(\alpha,\boldsymbol{\beta}),
\end{equation}
and
\begin{equation}\label{L_alpha_H}
\frac{\partial^2 \mathcal{L}_{\alpha}(\boldsymbol{\beta};\lambda ,\eta)}{\partial \boldsymbol{\beta}\partial \boldsymbol{\beta}^{T}} = \sum_{n=1}^N\boldsymbol{x}_n{\boldsymbol{x}_n}^{T} + \lambda(1-\eta) \boldsymbol{I}+ \lambda\eta \boldsymbol{H}(\alpha,\boldsymbol{\beta}),
\end{equation}
where $\boldsymbol{g}(\alpha,\boldsymbol{\beta})$ is the gradient of $\alpha$-absolute function which is defined elementwisely, i.e., $\boldsymbol{g}(\alpha,\boldsymbol{\beta})_{(j)} = \bigtriangledown |\beta_j|_\alpha$. $\boldsymbol{H}(\alpha,\boldsymbol{\beta})$ is a $p*p$ diagonal matrix and $\boldsymbol{H}(\alpha,\boldsymbol{\beta})_{(j,j)} = \bigtriangledown^2 |\beta_j|_{\alpha}, j=1,2,...,p$.
 \par
 It requires $O\left(\max(N,p)p^2\right)$ computing time in each iteration, then the full optimization procedure needs $O\left(T \max(N,p)p^2\right)$ computing time, where $T$ is the number of iterations required for the procedure to converge. Optimizer and statistical algorithm are blocked due to computational resource problems when the sample size $N$ becomes extremely large. In consequence, subsampling method is expected to exact useful and predominant information from massive data.

\section{Main Results: General Subsampling Algorithm and its Asymptotic Properties}
\label{sec:mainres}
In this section, we present a general subsampling algorithm for approximating the unknown optimal estimator $\widehat{\boldsymbol{\beta}}_{osa}$ and establish consistency and asymptotic normality of the resultant estimator. The general subsampling procedure is described in Algorithm~\ref{alg:general-alg}.

\begin{breakablealgorithm}
\caption{Randomized Smooth Elastic-net Regression}
\label{alg:general-alg}
{\bf Subsampling:} Form a subsample set $\{\widetilde{\boldsymbol{x}}_c,\widetilde{y}_c\}_{c = 1}^C$ of size $C$ by independently drawing from the original dataset with replacement according to sampling probabilities $\{\pi_n\}_{n=1}^N$.\\
{\bf Estimation:} Minimize the following penalized residual error function to get the subsampled estimator $\widetilde{\boldsymbol{\beta}}$ base on the drawn subsamples.\\
\begin{equation}
\widetilde{\mathcal{L}}_{\alpha}(\boldsymbol{\beta};\lambda ,\eta) =\frac{1}{2C}\sum_{c=1}^C{\frac{1}{\widetilde{\pi}_c}(\widetilde{\mathbf{x}}_c^\top\boldsymbol{\beta} - {\widetilde{y}_c})^2} + \frac{\lambda}{2} (1-\eta){\|\boldsymbol{\beta}\|}^2_2 + \lambda \eta{\|\boldsymbol{\beta}\|}_{\alpha}
\end{equation}
Due to the convexity of $\widetilde{\mathcal{L}}_{\alpha}(\boldsymbol{\beta};\lambda ,\eta)$ respect to $\boldsymbol{\beta}$, the minimization can be implemented by Newton's Method, iteratively applying the following formula until $\widetilde{\boldsymbol{\beta}}^{(t+1)}$ and $\widetilde{\boldsymbol{\beta}}^{(t)}$ are close enough.
\begin{equation}
\widetilde{\boldsymbol{\beta}}^{(t+1)} = \widetilde{\boldsymbol{\beta}}^{(t)} - \biggl[\frac{\partial^2 \widetilde{\mathcal{L}}_{\alpha}(\boldsymbol{\beta};\lambda ,\eta)}{\partial \boldsymbol{\beta} ^2}\big|_{\boldsymbol{\beta} = \widetilde{\boldsymbol{\beta}}^{(t)}}\biggl]^{-1}\frac{\partial \widetilde{\mathcal{L}}_{\alpha}(\boldsymbol{\beta};\lambda ,\eta,)}{\partial \boldsymbol{\beta}}\big|_{\boldsymbol{\beta} = \widetilde{\boldsymbol{\beta}}^{(t)}}
\end{equation}
\end{breakablealgorithm}
\par
Now we present two assumptions as the beginning of the analysis of consistency of $\widetilde{\boldsymbol{\beta}}$ to $\widehat{\boldsymbol{\beta}}_{osa}$ and asymptotic normality of $\widetilde{\boldsymbol{\beta}} - \widehat{\boldsymbol{\beta}}_{osa}$.
\begin{assumption}\label{secondMomentAss}
As $N \rightarrow{\infty}$, 
$$N^{-1}\sum_{n=1}^N{\|\boldsymbol{x}_n\|}^2 = O_P(1).$$
\end{assumption}
\par
\begin{assumption}\label{ISfourthMomentAss}
As $N \rightarrow{\infty}$,
$$N^{-2}\sum_{n=1}^N\pi_{n}^{-1}\|\boldsymbol{x}_n\|^4 = O_P(1).$$
\end{assumption}
\par
Assumption \ref{secondMomentAss} guarantees the existence of second moments of the full sample set. Assumption \ref{ISfourthMomentAss} pertains to both SSP and full sample set, in which $\frac{1/N}{\pi_n}$ can be regarded as the weighted ratio from the point of view of importance sampling (IS) \citep{neal2001annealed}.

\begin{thm}\label{thm2}
If Assumption~\ref{secondMomentAss} and~\ref{ISfourthMomentAss} are satisfied, as $N\rightarrow{\infty}$ and $C\rightarrow{\infty}$, $\widetilde{\boldsymbol{\beta}}$ is consistent to $\widehat{\boldsymbol{\beta}}_{osa}$ in conditional probability, given $\mathcal{F}_N$ in probability. Moreover, the rate of convergence is $C^{-1/2}$. That is, with probability approaching one, for any $\epsilon > 0$, there exists $\Delta_{\epsilon}$ and $C_{\epsilon}$, such that for all $C>C_{\epsilon}$\\
\begin{equation}
\label{3.15}
P\left(\|\widetilde{\boldsymbol{\beta}} - \widehat{\boldsymbol{\beta}}_{osa}\|\geq C^{-1/2}\Delta_{\epsilon}|\mathcal{F}_N\right) < \epsilon.
\end{equation}

\end{thm}
\par
Theorem \ref{thm2} demonstrates that $\widetilde{\boldsymbol{\beta}}$ converges to $\widehat{\boldsymbol{\beta}}_{osa}$ in probability as $C\rightarrow\infty$ and $N\rightarrow\infty$. Additionally, Theorem \ref{thm1} and Theorem \ref{thm2} indicate that $\widetilde{\boldsymbol{\beta}}$ is a biased estimate of $\widehat{\boldsymbol{\beta}}_{opt}$ and the bias diminishes to zero as $\alpha\rightarrow\infty$. Hence the analysis of asymptotic properties based on smooth Elastic-net regression is of referential significance for the original Elastic-net regression in condition of large $\alpha$.

\begin{assumption}\label{ISdeltaMomentAss}
There exists some $\delta > 0$, such that 
$$N^{-(2+\delta)}\sum_{n=1}^N{\pi_n}^{-1-\delta}\|\boldsymbol{x}_n\|^{2\delta + 4} = O_P(1).$$
\end{assumption}
For convenience, denote 
\begin{equation}\label{MxDef}
\mathbf{M_X} = \frac{1}{N}\frac{\partial^2 \mathcal{L}_{\alpha}(\boldsymbol{\beta};\lambda ,\eta)}{\partial \boldsymbol{\beta} ^2}|_{\boldsymbol{\beta}=\widehat{\boldsymbol{\beta}}_{osa}},
\end{equation}
\begin{equation}
\mathbf{C}_{osa} = \lambda(1-\eta)\widehat{\boldsymbol{\beta}}_{osa} + \lambda\eta\boldsymbol{g}(\alpha, \widehat{\boldsymbol{\beta}}_{osa}).
\end{equation}
If $\lambda>0$ and $\eta \in (0,1)$,it is obvious that $\mathbf{M_X}$ is a positive definite matrix from \eqref{alphaAbsH} and \eqref{L_alpha_H}. Moreover, according to first-order necessary condition from optimization theory and \eqref{L_alpha_G},
\begin{equation}\label{fordc}
\frac{\partial \mathcal{L}_{\alpha}(\boldsymbol{\beta};\lambda ,\eta)}{\partial \boldsymbol{\beta}}|_{\boldsymbol{\beta}=\widehat{\boldsymbol{\beta}}_{osa}} = \sum_{n=1}^N(\widehat{\boldsymbol{\beta}}_{osa}^{T} \boldsymbol{x}_n - y_n)\boldsymbol{x}_n + \mathbf{C}_{osa} = 0.  
\end{equation}
The below theorem presents the asymptotic normality of $\widetilde{\boldsymbol{\beta}} - \widehat{\boldsymbol{\beta}}_{osa}$.

\begin{thm}\label{ANoftildeBeta}
If Assumption~\ref{secondMomentAss},~\ref{ISfourthMomentAss} and~\ref{ISdeltaMomentAss} hold, then as  $N \rightarrow{\infty}$ and  $C\rightarrow{\infty}$, conditional on $\mathcal{F}_N$ in probability, 
\begin{equation}
\label{AN}
{\bf V}^{-1/2}(\widetilde{\boldsymbol{\beta}} - \widehat{\boldsymbol{\beta}}_{osa}) \rightarrow\mathcal{N}\left(0, \mathbf{I}\right)
\end{equation}
in distribution, where
\begin{equation}
\label{Vform}
{\bf V} = \mathbf{M_X}^{-1}\mathbf{V}_0\mathbf{M_X}^{-1} = O_{P|\mathcal{F}_N}(C^{-1})
\end{equation}
and
\begin{equation}\label{V0form}
\begin{aligned}
\mathbf{V}_0 = \frac{1}{CN^2}\biggl[ \sum_{n=1}^N&{\pi_n}^{-1}({\widehat{\boldsymbol{\beta}}_{osa}^\top\boldsymbol{x}_n-y_n})^2\boldsymbol{x}_n{\boldsymbol{x}_n}^{\top} + \mathbf{C}_{osa}{\mathbf{C}_{osa}}^{\top} \\
&+ \sum_{n=1}^N(\widehat{\boldsymbol{\beta}}_{osa}^\top\boldsymbol{x}_n-y_n)(\mathbf{C}_{osa}{\boldsymbol{x}_n}^{\top}+\boldsymbol{x}_n\mathbf{C}_{osa}^{\top}) \biggl].
\end{aligned}
\end{equation}
\end{thm}

\section{Optimal Subsampling Probabilities Design}\label{sec:optimal}
In accordance with Theorem \ref{ANoftildeBeta}, the asymptotic variance of $\widetilde{\boldsymbol{\beta}}$ obtained from general subsampling algorithm \ref{alg:general-alg} is derived. Then we aim to design OSP by minimizing the asymptotic MSE of $\widetilde{\boldsymbol{\beta}}$.

\begin{thm}[Optimal Subsampling Probability]
\label{thm4}
In Algorithm~\ref{alg:general-alg}, if SSP is chosen such that $\pi_n^{osp}\propto\pi_n^*$ with $\pi_n^*$ satisfying the following quadratic equation 
\begin{equation*}
\pi^*_n = \|\mathbf{M_X^{-1}}\mathbf{Z}_n\|^2, \text{ and }
\mathbf{Z}_n =  \pi_n^{*-1/2}(\boldsymbol{x}_n^{\top}{\widehat{\boldsymbol{\beta}}_{osa}}-y_n)\boldsymbol{x}_n + \pi_n^{*1/2}\mathbf{C}_{osa},
\end{equation*}
for $n = 1,2,\dots,N$, the asymptotic MSE of $\widetilde{\boldsymbol{\beta}}_{osa}$, i.e., $\text{tr}(\mathbf{V})$ attains its minimum. 
\end{thm}
\par
The detailed explanation of the optimal subsampling condition in Theorem \ref{thm4} is that $\{\pi_n^{osa}\}_{n=1}^N$ satisfy the following equations:
\begin{equation}
\begin{split}
\label{QEs}
(\mathbf{M_X}^{-1}\mathbf{C}_{osa})^{T}(\mathbf{M_X}^{-1}\mathbf{C}_{osa})(\pi_n^{osp})^{2} + 2(\widehat{\boldsymbol{\beta}}_{osa}^{T}\boldsymbol{x}_n - y_n)(\mathbf{M_X}^{-1}\mathbf{C}_{osa})^{T}(\mathbf{M_X}^{-1}\boldsymbol{x}_n)\pi_n^{osp} \\
+ (\widehat{\boldsymbol{\beta}}_{osa}^{T}\boldsymbol{x}_n - y_n)^2(\mathbf{M_X}^{-1}\boldsymbol{x}_n)^{T}(\mathbf{M_X}^{-1}\boldsymbol{x}_n) = \mathcal{K}(\pi_n^{osp})^{2},~~~n=1,2,...,N,
\end{split}
\end{equation}
where $\mathcal{K}$ is a proportionality coefficient. The other hidden conditions that $\{\pi_n^{osp}\}_{n=1}^N$ have to satisfy are
\begin{equation}
\label{sumE}
\sum_{n=1}^N\pi_n^{osp} = 1
\end{equation}
and
\begin{equation}\label{nonN}
\pi_n^{osp} \geq 0, n=1,2,...,N.
\end{equation}
\par
According to \eqref{QEs}, \eqref{sumE} and \eqref{nonN}, the problem of searching OSP for smooth Elastic-net is translated to seeking a suitable proportionality coefficient $\mathcal{K}$ which ensures that all quadratic equations in \eqref{QEs} have at least one non-negative solution and the set composed of the N non-negative solutions satisfies \eqref{sumE} and \eqref{nonN}.
\par
Denote the quadratic coefficient in \eqref{QEs} as
\begin{equation}\label{2co}
\mathcal{A} = (\mathbf{M_X}^{-1}\mathbf{C}_{osa})^{T}(\mathbf{M_X}^{-1}\mathbf{C}_{osa}) - \mathcal{K}.
\end{equation}
Denote the primary coefficients in \eqref{QEs} as
\begin{equation}\label{1co}
\mathcal{B}_n = 2(\widehat{\boldsymbol{\beta}}_{osa}^{T}\boldsymbol{x}_n - y_n)(\mathbf{M_X}^{-1}\mathbf{C}_{osa})^{T}(\mathbf{M_X}^{-1}\boldsymbol{x}_n).
\end{equation}
Denote the constant terms in \eqref{QEs} as
\begin{equation}\label{cco}
\mathcal{D}_n = (\widehat{\boldsymbol{\beta}}_{osa}^{T}\boldsymbol{x}_n - y_n)^2(\mathbf{M_X}^{-1}\boldsymbol{x}_n)^{T}(\mathbf{M_X}^{-1}\boldsymbol{x}_n).
\end{equation}
\par
Based on the formulas above, we get the property that $\mathcal{D}_n \geq 0 ,~n =1,2,...,N.$ When $\mathcal{A} < 0$, i,e., $\mathcal{K} >(\mathbf{M_X}^{-1}\mathbf{C}_{osa})^{T}(\mathbf{M_X}^{-1}\mathbf{C}_{osa})$, every quadratic equation in \eqref{QEs} has at least one non-negative solution, which can be expressed in the following form:
 \begin{equation}\label{wdpi}
 \pi_n^{osp} = \frac{-\mathcal{B}_n + \sqrt{\mathcal{B}_n^2 - 4\mathcal{A}\mathcal{D}_n}}{2\mathcal{A}} = -\frac{\mathcal{B}_n}{\mathcal{A}} + \sqrt{\frac{1}{4}(\frac{\mathcal{B}_n}{\mathcal{A}})^2 - \frac{\mathcal{D}_n}{\mathcal{A}}}, ~n = 1,2,...,N.
 \end{equation}
 
Since $\{\mathcal{B}_n\}$ and $\{\mathcal{D}_n\}$ depend on full data set, even if we find a proportional coefficient $\mathcal{K}$ that ensures the existence of at least one non-negative solution to all quadratic equations in \eqref{QEs}, we still cannot guarantee that these non-negative solutions satisfy \eqref{sumE}, i.e., the sum equals 1.
So it is hard to find such $\mathcal{K}$
that simultaneously satisfies \eqref{QEs}, \eqref{sumE}, and \eqref{nonN} (even such $\mathcal{K}$ does not exist). In consequence, it's necessary to find a suitable substitution for \eqref{QEs} (or \eqref{wdpi}) under a feasible $\mathcal{K}$.
\\
\par
\textbf{Pseudo-optimal SSP.}
Consider the proportionality coefficient as 
\begin{equation}\label{chosenK}
\mathcal{K} = (\mathbf{M_X}^{-1}\mathbf{C}_{osa})^{T}(\mathbf{M_X}^{-1}\mathbf{C}_{osa}) + (\sum_{n=1}^N |\widehat{\boldsymbol{\beta}}_{osa}^{T}\boldsymbol{x}_n - y_n|~||\mathbf{M_X}^{-1}\boldsymbol{x}_n||_2)^2. 
\end{equation}
Then we have
$$\mathcal{A} = -(\sum_{n=1}^N |\widehat{\boldsymbol{\beta}}_{osa}^{T}\boldsymbol{x}_n - y_n|~||\mathbf{M_X}^{-1}\boldsymbol{x}_n||_2)^2 \leq 0.$$
Get back to \eqref{wdpi}, it is observed that if the term $\mathcal{B}_n$ is neglected, then we obtain a substitution for \eqref{wdpi} which simultaneously satisfies condition \eqref{sumE} and condition \eqref{nonN}, i.e., 
\begin{equation}
\label{apposp}
\pi_n^{osp}\approx\pi_n^{posp}=\sqrt{-\frac{\mathcal{D}_n}{\mathcal{A}}} = \frac{|\widehat{\boldsymbol{\beta}}_{osa}^{T}\boldsymbol{x}_n - y_n|~||\mathbf{M_X}^{-1}\boldsymbol{x}_n||_2}{\sum_{n=1}^N |\widehat{\boldsymbol{\beta}}_{osa}^{T}\boldsymbol{x}_n - y_n|~||\mathbf{M_X}^{-1}\boldsymbol{x}_n||_2}, ~~~~n=1,2,...,N.
\end{equation}
We denote the SSP $\boldsymbol{\pi}^{posp}$ in \eqref{apposp} as \textbf{pseudo-optimal subsampling probability (POSP)} for smooth Elastic-net. 
\par
In Algorithm \ref{alg:general-alg}, if uniform probability is chosen as SSP, we represent the corresponding asymptotic MSE as $tr(\mathbf{V}^{uni})$. If $\boldsymbol{\pi}^{posp}$ is chosen as SSP, we represent the corresponding asymptotic MSE as $tr(\mathbf{V}^{posp})$. 
The following corollary demonstrates that, although $\boldsymbol{\pi}^{posp}$ cannot achieve the same performance as $\boldsymbol{\pi}^{osp}$ in Theorem \ref{thm1}, it still outperforms uniform SSP in terms of MSE.
\begin{corollary}\label{coro1}
$tr(\mathbf{V}^{uni}) \geq tr(\mathbf{V}^{posp}).$
\end{corollary}
\par
From the proof in supplementary material, it is also observed that under the condition where the distribution of $\boldsymbol{x}$ is relatively dispersed (i.e., the heavy-tailed distribution), the superiority of $\{\pi_n^{posp}\}$ over uniform SSP will be more significant.
\section{Implementation}\label{sec:twostep}
To implement the Randomized Smooth Elastic-net Regression Algorithm~\ref{alg:general-alg} with $\boldsymbol{\pi}^{posp}$ in \eqref{apposp}, the following two-step algorithm is developed. 

\begin{breakablealgorithm}
\caption{Two-Step Algorithm}
\label{alg:2step-alg}
{\bf Input}: Subsample size $C_0$ and $C$. \\
{\bf Step 1}: Run Algorithm~\ref{alg:general-alg} with subsample size $C_0$ and uniform sampling probability $\left\{\pi_n = \frac{1}{N}\right\}_{n=1}^N$ to obtain the pilot estimate of $\widehat{\boldsymbol{\beta}}_{osa}$ denoted by $\widetilde{\boldsymbol{\beta}}_0$. Replace $\widehat{\boldsymbol{\beta}}_{osa}$ with $\widetilde{\boldsymbol{\beta}}_0$ in \eqref{MxDef} and \eqref{apposp} to calculate the nearly exact approximation of optimal SSP.\\
{\bf Step 2}: Run Algorithm~\ref{alg:general-alg} with subsample size $C$ and the SSP obtained from step 1 to obtain the estimator $\breve{\boldsymbol{\beta}}$. \\
\end{breakablealgorithm}

\section*{Asymptotic Properties}
Now we consider the consistency of $\breve{\boldsymbol{\beta}}$ based on $\{\pi_n^{posp}(\widetilde{\boldsymbol{\beta}}_0)\}_{n=1}^N$ to $\widehat{\boldsymbol{\beta}}_{osa}$. We establish its asymptotic properties with the following assumption. Note that all subsampling probabilities of $\widetilde{\mathcal{L}_{\alpha}}$ in following statement is replaced with $\{\pi_n^{posp}(\widetilde{\boldsymbol{\beta}}_0)\}_{n=1}^N$. To simplify the symbol, we will not make special annotations in the following formulas.
\begin{assumption}\label{ass4}
The covariate distribution satisfies that $E(\|\boldsymbol{x}\|^4) < \infty$.  
\end{assumption}
Assumption \ref{ass4} guarantees that the second and fourth moments of the covariate distribution exist.
\begin{thm}\label{brevebetaconsistency}
Under Assumption \ref{ass4}, if the estimate $\widetilde{\boldsymbol{\beta}}_0$ based on the first step sample exits, then as $C\rightarrow\infty$ and as $N\rightarrow\infty$, with probability approaching one, for any $\epsilon > 0$, there exists a finite $\Delta_{\epsilon}$ and $C_{\epsilon}$, such that for all $C > C_{\epsilon}$,
\begin{equation} 
P\left(\|\breve{\boldsymbol{\beta}} - \widehat{\boldsymbol{\beta}}_{osa}\|\geq C^{-1/2}\Delta_{\epsilon}|\mathcal{F}_N\right) < \epsilon.
\end{equation}

\end{thm}

\begin{thm}\label{brevebetaAN}
Under Assumption \ref{ass4}, as $C_0 \rightarrow \infty$ , $C\rightarrow \infty$ and $N\rightarrow\infty$, conditional on $\mathcal{F}_N$ and $\widetilde{\boldsymbol{\beta}}_0$,
\begin{equation}\label{breveBetaAN}
\mathbf{V}^{\widehat{\boldsymbol{\beta}}_{osa}-1/2}(\breve{\boldsymbol{\beta}} - \widehat{\boldsymbol{\beta}}_{osa})\rightarrow\mathcal{N}(0, \mathbf{I})
\end{equation}
in distribution, where $\mathbf{V}^{\widehat{\boldsymbol{\beta}}_{osa}} = \mathbf{M_X}^{-1}\mathbf{V}_0^{\widehat{\boldsymbol{\beta}}_{osa}}\mathbf{M_X}^{-1}$ and 
\begin{equation}
\begin{split}
\mathbf{V}_0^{\widehat{\boldsymbol{\beta}}_{osa}} = \frac{1}{CN^2}\biggl[ \sum_{n=1}^N|\widehat{\boldsymbol{\beta}}_{osa}^{\top}\boldsymbol{x}_n - y_n|\|\mathbf{M_X}^{-1}\boldsymbol{x}_n\|\sum_{n=1}^N\frac{|{\widehat{\boldsymbol{\beta}}_{osa}^{\top}\boldsymbol{x}_n-y_n}|}{\|\mathbf{M_X}^{-1}\boldsymbol{x}_n\|}\boldsymbol{x}_n{\boldsymbol{x}_n}^{\top}\\ + \mathbf{C}_{osa}{\mathbf{C}_{osa}}^{\top}
+ \sum_{n=1}^N(\widehat{\boldsymbol{\beta}}_{osa}^{\top}\boldsymbol{x}_n-y_n)(\mathbf{C}_{osa}{\boldsymbol{x}_n}^{\top}+\boldsymbol{x}_n\mathbf{C}_{osa}^{\top}) \biggl]
\end{split}
\end{equation}
\end{thm}
\par
It is worth noting that Theorem \ref{brevebetaconsistency} and Theorem \ref{brevebetaAN} hold under Assumption \ref{ass4}, which is much weaker than the prerequisite conditions of Theorem \ref{thm2} and Theorem \ref{ANoftildeBeta}. This implies that utilizing the subsampling probabilities in the form of (\ref{apposp}) promotes the attainment of consistency and asymptotic normality of the subsampling estimator. In practical applications, to avoid large estimation error of POSP in step 1 such that severe estimation bias occurs on the estimator $\breve{\boldsymbol{\beta}}$ in step 2, $C_0$ is not advisable to be too small. The issue of $C_0$ proportion will be further investigated in \S \ref{sec:simulation}.

\section{Simulation}\label{sec:simulation}
The main purpose of this simulation study is to verify the theoretical results of Algorithm \ref{alg:2step-alg} proposed in \S\ref{sec:twostep}.
The data which we simulate is generated from the true model:
$$\mathbf{Y} = \mathbf{X}\boldsymbol{\beta} + \sigma\mathbf{\epsilon},$$
where
$$\mathbf{\epsilon}\sim\mathcal{N}(0, \mathbf{I}),$$
and careful selection of $\sigma$ is necessary depending on different cases.
\par
Within each case, we let $\alpha = 10$ and the number of predictors $p =50$. We use 500 observations to conduct K-fold cross-validation ($K = 5$) for hyper-parameters $\lambda$ and $\eta$ on a two-dimension surface in the shape of $\{e^{-3},e^{-2},e^{-1},e^{0},e^{1},...,e^{10}\} * \{0.1,0.2,0.3,...,0.9\}$. Considering the positive definiteness of $\mathbf{M_X}$ and strict convexity of $\mathcal{L}_{\alpha}$, the hyper-parameters must meet the following conditions: $\lambda \neq 0$ and $\eta \neq 0$.
\par
Here are details of four examples:
\begin{itemize}
    \item \textbf{case 1:} $\boldsymbol{x}\sim\mathcal{N}(0,\mathbf{\Sigma})$, $\mathbf{\Sigma}^{(i,j)} = 0.5^{|i-j|}$, \\$\boldsymbol{\beta} = (\overbrace{4,...,4}^{10},\overbrace{0,...,0}^{10},\overbrace{2,...,2}^{10},\overbrace{0,...0}^{10},\overbrace{1,...,1}^{10})$ and $\sigma = 3$. Case 1 reproduces the most common application scenario of Elastic-net, in which the target coefficients exhibit sparsity and linear model is independent of some features. Under such data distribution, Assumption \ref{ass4} holds naturally.
    \item \textbf{case 2:} Case 2 is the same as case 1 except that $\boldsymbol{x}\sim\mathbf{t}_3(0,\mathbf{\Sigma})$. For this case, the distribution of $\boldsymbol{x}$ has a heavy tail and does not satisfy Assumption \ref{ass4}. This case is designed to examine how sensitive Algorithm \ref{alg:2step-alg} is to Assumption \ref{ass4}.
    \item \textbf{case 3:} $\boldsymbol{x}_i$ are independent identically distributed $exp(2.0),i=1,...,50$, $\boldsymbol{\beta} = (2,...,2)$ and $\sigma = 3$. For this case, the distribution of $\boldsymbol{x}$ has a heavier tail on the right and all features are equally important.
    \item \textbf{case 4:} 
    $$\boldsymbol{x}_i = Z_1 + \epsilon_i^{x},  Z_1\sim\mathcal{N}(0, 1),i = 1,2,3,4,5,$$
    $$\boldsymbol{x}_i = Z_2 + \epsilon_i^{x},  Z_2\sim\mathcal{N}(0, 1),i = 6,7,8,9,10,$$
    $$\boldsymbol{x}_i = Z_3 + \epsilon_i^{x},  Z_3\sim\mathcal{N}(0, 1),i = 11,12,13,14,15,$$
    $$\boldsymbol{x}_i = Z_4 + \epsilon_i^{x},  Z_4\sim\mathcal{N}(0, 1),i = 16,17,18,19,20,$$
    where $\epsilon_i^{x}$ are independent identically distributed $\mathcal{N}(0,0.01), i =1,...,20$. 
    $\boldsymbol{x}_i\sim\mathcal{N}(0,1),$ $\boldsymbol{x}_i$ is independent identically distributed, $i = 21,...,50.$
    $\boldsymbol{\beta} = (\overbrace{3,...3}^{20},\overbrace{0,...,0}^{30})$ and $\sigma = 3$.
    In this case, we have four equally important groups composed of five features and 30 noise features.
\end{itemize}
\par
To sum up, data distributions of case 1 and case 4 have light tails while data distributions of case 2 and case 3 are heavy-tailed. Note that case 4 is a grouped variable situation designed by \citet{zou2005regularization} to examine the coefficients path of Elastic-net regression.
\par
\begin{figure}[htbp]
    \centering
    \includegraphics[width=1.0\columnwidth]{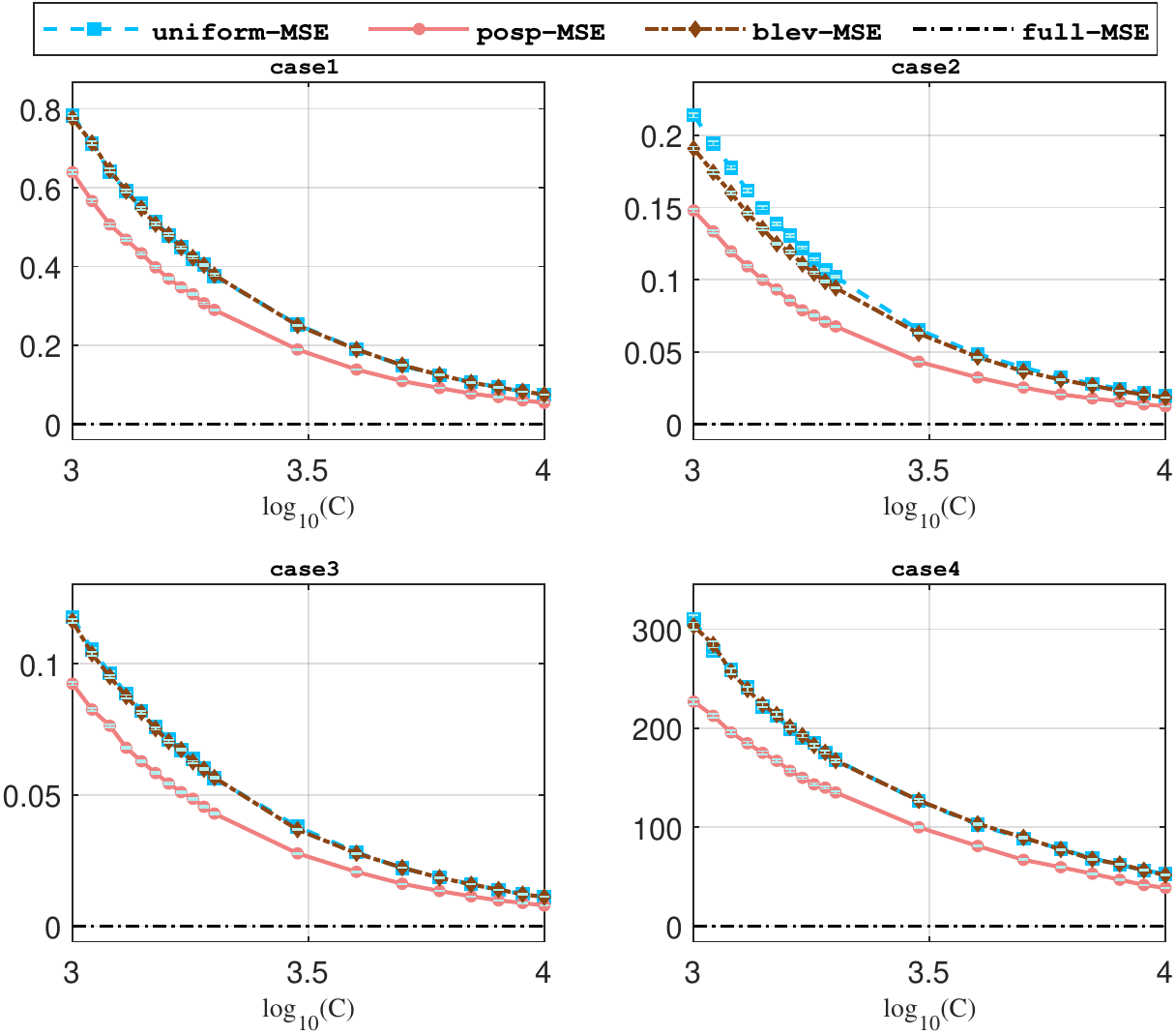}
    \caption{MSEs for different $C$ for step 2 in Algorithm \ref{alg:2step-alg}. The subsample size $C_0$ for step 1 is fixed to 1000. The subsample size for uniform and BLEV is the same as $C$.}
    \label{MSE1_fig}
\end{figure}

To evaluate the consistency of $\breve{\boldsymbol{\beta}}$ to full data set regression estimator $\widehat{\boldsymbol{\beta}}_{osa}$, for Algorithm \ref{alg:2step-alg}, subsample size $C_0$ for step 1 will be fixed to 1000. Subsample size $C$ for step 2 is going to increase from 1000 to 1900 with step size 100 and increase from 2000 to 10000 with step 1000.  Thereby the standard optimal estimator $\widehat{\boldsymbol{\beta}}_{osa}$ will be established via running Newton's Method on the full data set, which has the scale of $N = 1e5$. For comparison, uniform subsampling and BLEV will also be applied. BLEV was proposed by \citet{drineas2011faster}. If $\mathbf{X}$ is a $N*p$ sample matrix and $\mathbf{U}$ is a $N*p$ orthogonal matrix obtained by SVD decomposition for $\mathbf{X}$, whose columns contain the left singular vectors of $\mathbf{X}$. Assume that $\boldsymbol{u}_{n}$ is the $\it{n}$-th row of $\mathbf{U}$, then BLEV probability has the following form:
$$\pi_{n}^{blev} = \frac{\|\boldsymbol{u}_{n}\|^2}{\sum_{n=1}^N\|\boldsymbol{u}_{n}\|^2}, ~~~n = 1,...,N.$$
Experiment under each $C$ will be repeated for $M = 1000$ times. Mean squared error (MSE) which measures subsampling consistency will be calculated:
$$MSE = \frac{1}{M}\sum_{k=1}^{M}\|\breve{\boldsymbol{\beta}}^{(k)} - \widehat{\boldsymbol{\beta}}_{osa}\|^2.$$
MSE of uniform SSP and BLEV SSP will be calculated in the same way. The result is shown in Fig \ref{MSE1_fig}. For all four cases, SSP $\boldsymbol{\pi}^{posp}$ in form of \eqref{apposp} always performs better than $\boldsymbol{\pi}^{uni}$ and $\boldsymbol{\pi}^{blev}$, which agrees with the theory of Corollary \ref{coro1}. In smooth Elastic-net regression scenario, the commonly used $\boldsymbol{\pi} ^{blev}$ does not exhibit significant advantage over uniform subsampling.
\par
\begin{figure}[htbp]
    \centering
    \includegraphics[width=1.0\columnwidth]{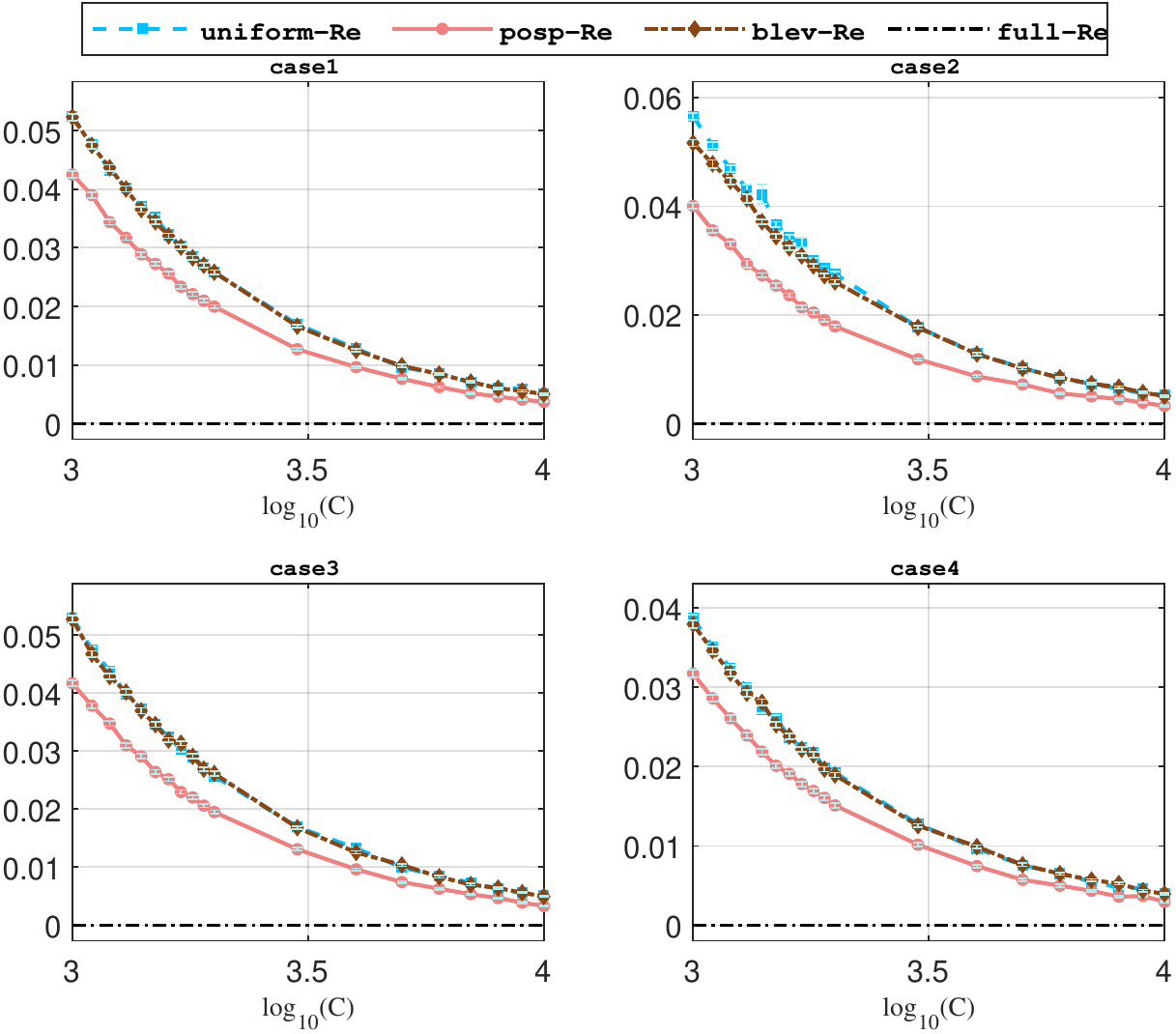}
    \caption{Res for different $C$ for step 2 in Algorithm \ref{alg:2step-alg}. The subsample size $C_0$ for step 1 is fixed to 1000. The subsample size for uniform and BLEV is the same as $C$.}
    \label{Re1_fig}
\end{figure}

Furthermore, we are concerned with the relative-error approximation Re:
$$Re = \frac{1}{M}\sum_{k=1}^{M}\frac{\|\mathbf{X}_{test}\breve{\boldsymbol{\beta}}^{(k)} - \mathbf{Y}_{test}\|^2}{\|\mathbf{X}_{test}\boldsymbol{\widehat{\boldsymbol{\beta}}}_{osa} - \mathbf{Y}_{test}\|^2} - 1,$$
where the scale of test data set $\{\mathbf{X}_{test}, \mathbf{Y}_{test}\}$ is 1000. The $Re$ result is shown in Fig \ref{Re1_fig}, which has the same trend as the $MSE$ result in Fig \ref{MSE1_fig}.
\par
\begin{figure}[htbp]
    \centering
    \includegraphics[width=1.0\columnwidth]{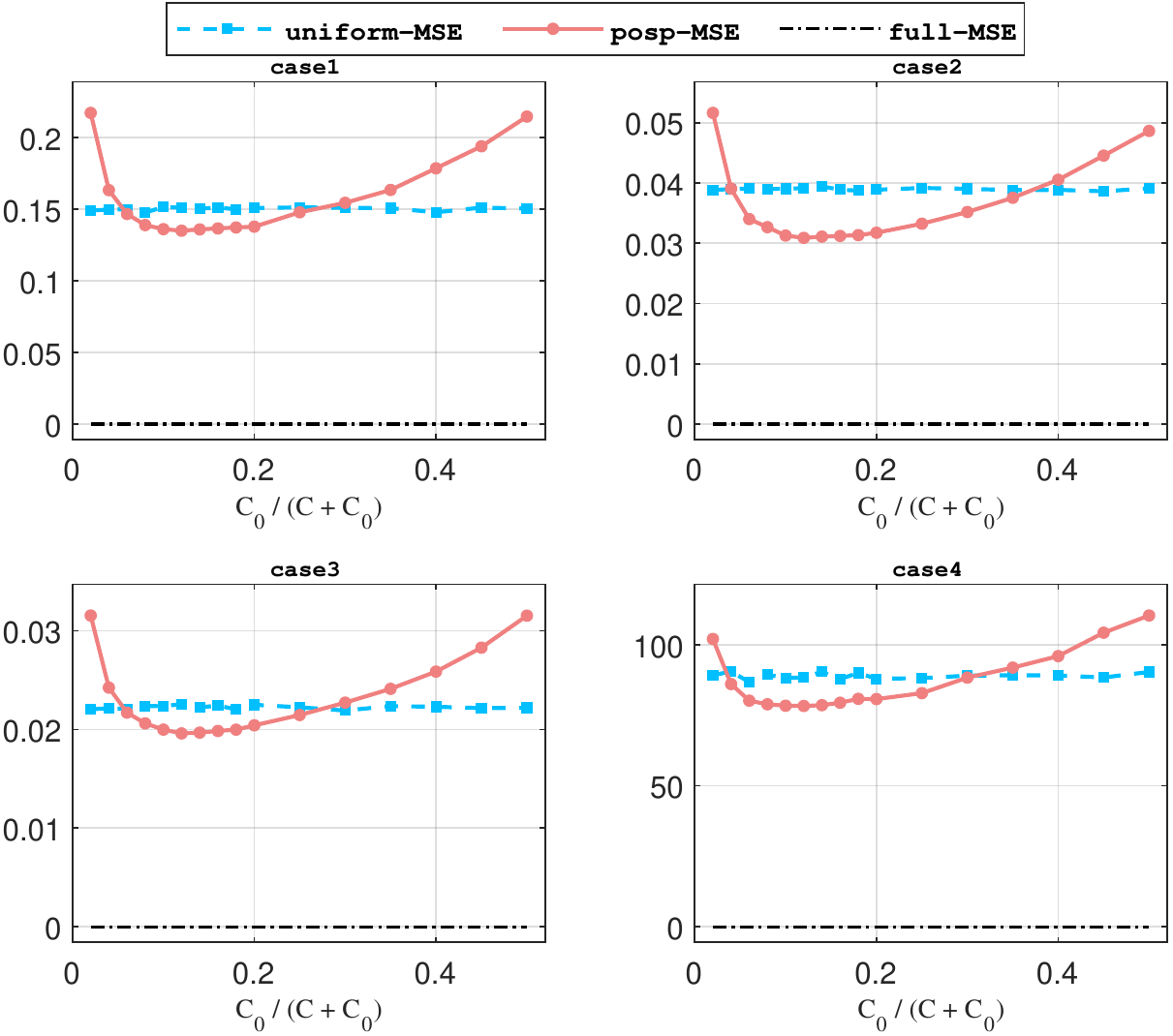}
    \caption{MSEs for different $C_0 / (C+C_0)$ for Algorithm \ref{alg:2step-alg}. $C + C_0$ is fixed to 5000. The subsample size for uniform is fixed to $C + C_0 = 5000$.}
    \label{MSE2_fig}
\end{figure}

To investigate the impact of sample allocation between step 1 and step 2 of Algorithm \ref{alg:2step-alg} on the final estimator, we fixed the total sample size at 5000, i.e., $C + C_0 = 5000.$ Let the proportion sequence of $C_0$ be $\{0.02,0.04,0.06,...,0.18,0.2,0.3,0.4,0.5\}$. As reference, we fix the subsample number of uniform subsampling to 5000. Following the same experimental setup as before, the performance of Algorithm \ref{alg:2step-alg} was observed under different proportions of $C_0$. As shown in Fig \ref{MSE2_fig}, Algorithm \ref{alg:2step-alg} performs the best when $C_0 / (C + C_0)$ is around 0.15. However, this is only the analysis result for the above four cases. When the total number of samples is limited, allocation ratio of $C_0$ under other scenarios still needs further investigation.
\par
All simulation experiments are carried out by matlab software on a computer equipped with 2.30 GHz intel i7 processor, 16.0 GB memory and Windows 11 operation system.

\section{Real Data Experiment}\label{sec:realdata}
\subsection{Blog feedback prediction}
In this section, we apply Algorithm \ref{alg:2step-alg} to BlogFeedback data set (\url{https://archive.ics.uci.edu/ml/datasets/BlogFeedback}) from the UCI repository\citep{Dua:2019}. The goal is to predict the feedback number of a given blog with 280 predictors including basic features (number and decrease/increase speed of links and feedbacks in past), textual features (the most discriminative bag of words features), weekday features (binary features that describe on which day of a week the blog was published) and parent features (parent blogs and the number of feedbacks they received). More details can be found in \citet{buza2013feedback}. 
\par
\begin{figure}[htbp]
    \centering
    \includegraphics[width=1.0\columnwidth]{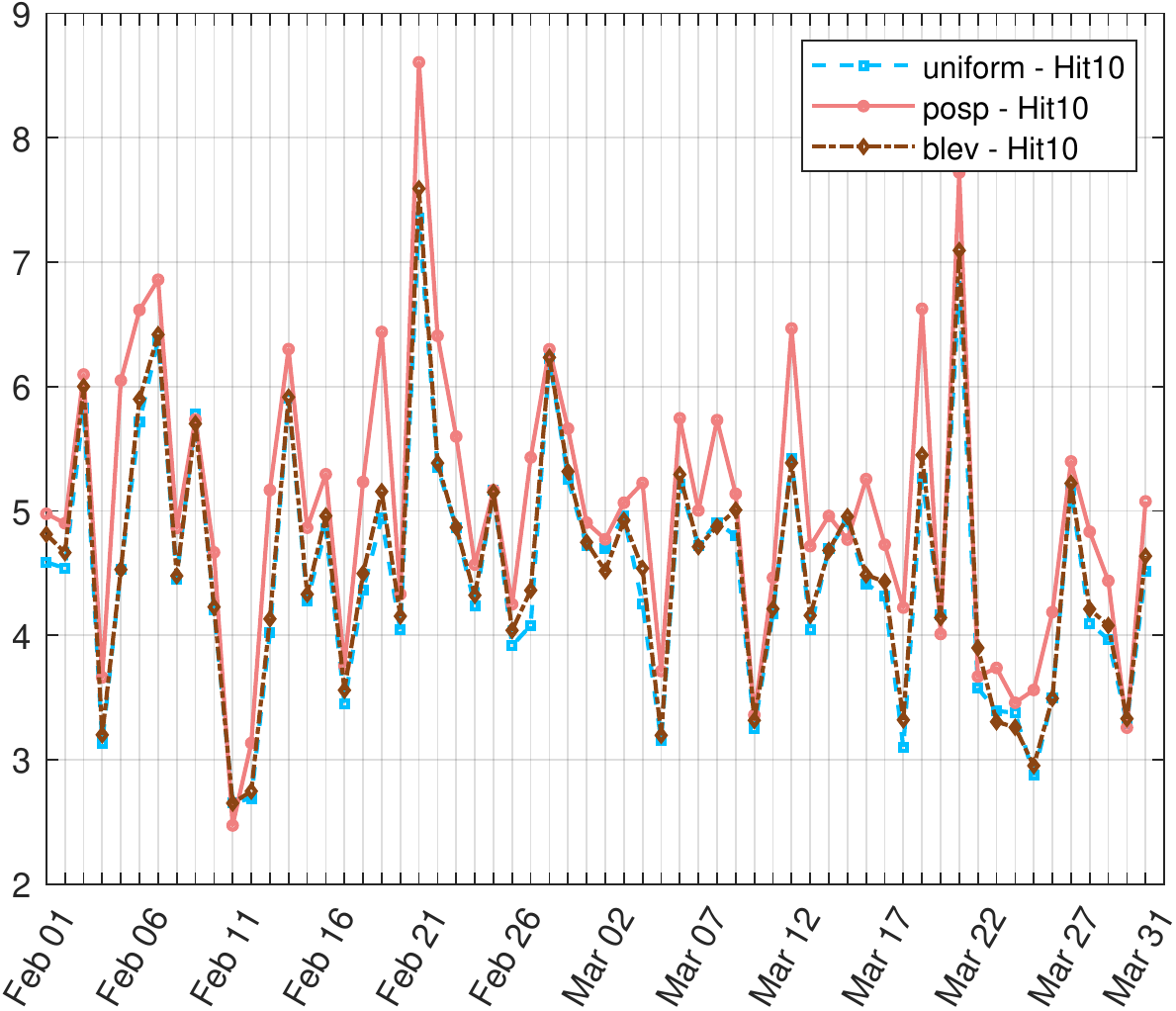}
    \caption{Average Hit-10 curves on test data set. Full train data set is from 2010 and 2011 and test data set is from 2012 Feb 1 to 2012 Mar 31.}
    \label{Hit10_fig}
\end{figure}

There are a total of 52397 train samples, which represent blog data from the years 2010 to 2011. We set $C_0 = 1000$ for step 1 and $C = 2000$ for step 2 in Algorithm \ref{alg:2step-alg}, and the resultant estimator will be tested on 60-day test data set from 2012 Feb 1 to 2012 Mar 31. Moreover, hyper-parameters $\lambda$, $\eta$ and $\alpha$ are fixed to $e^8$, 0.8 and 10. For each day in test set, 10 blogs will be predicted to have the largest number of feedbacks, then we count how many blogs out of these that recieved the largest number in reality. Such evaluation measure is called Hit-10. We repeat subsampling algorithm for 1000 times. Fig \ref{Hit10_fig} shows the plot of average Hit-10 of uniform, BLEV and POSP SSPs. Hit-10 curve of POSP is always above those of uniform and BLEV, indicating that POSP extracts more information from full data set and gives a more accurate prediction when subsample size is fixed.
\subsection{2D CT slice axial axis location prediction}
\begin{figure}[htbp]
    \centering
    \subfigure[Scaled CT image.]{\includegraphics[width=0.31\textwidth]{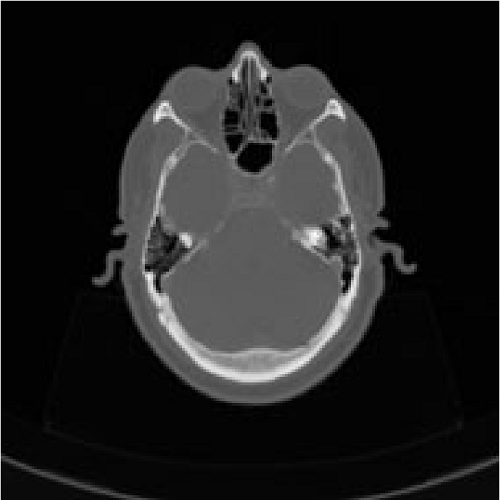} \label{subfig1:CT}}
    \subfigure[Detected region.]{\includegraphics[width=0.31\textwidth]{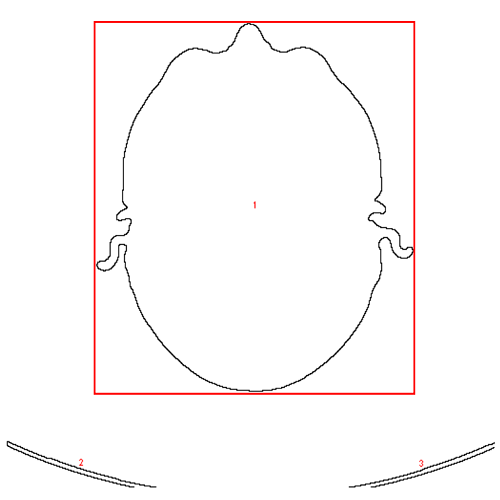} \label{subfig2:CT}}
    \subfigure[Sector/shell model.]{\includegraphics[width=0.31\textwidth]{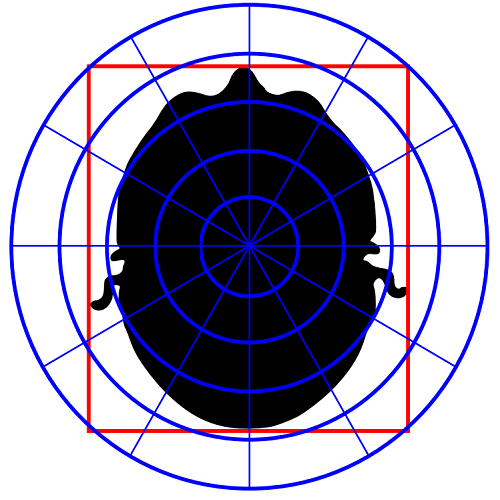} \label{subfig3:CT}}
    \caption{The feature extraction process of a neck scan image.\citep{graf20112d}}
    \label{fig:CT}
\end{figure}

In this section, we apply Algorithm \ref{alg:2step-alg} to the CT slices data set(\url{https://archive.ics.uci.edu/ml/datasets/Relative+location+of+CT+slices+on+axial+axis}) from the UCI repository\citep{Dua:2019}. The data was retrieved from a set of 53500 CT images from 74 different patients (43 male, 31 female). There are 384 predictors totally, which were exacted from the processed images by three steps as shown in Fig \ref{fig:CT}. First, it is needed to scale the CT image to a common resolution (1.5 px/mm), as shown in Fig \ref{subfig1:CT}. Second, a compound region detection algorithm was applied to separate the body region from the CT image, as shown in Fig \ref{subfig2:CT}. Finally, a radial sector/shell model was created from which the two descriptors representing dense structures (bones) and soft tissues (like organs etc.) are extracted, as shown in Fig \ref{subfig3:CT}. More details can be found in \citet{graf20112d}.
\par
The target is to predict the axial axis location between 0 and 180 (where 0 represents the top of head and 180 represents the soles of feet) of a CT slice to help radiologists localize and align parts of CT volume scans to perform among others differential diagnosis. We apply Algorithm \ref{alg:2step-alg} to predict relative location of a CT slice. We set $C_0 = 2000$ and let $C$ be 100, 1000 and 10000. As reference, Algorithm \ref{alg:general-alg} equipped with uniform and BLEV SSPs will also be applied under subsample size $C = 100, 1000, 10000$. Hyper-parameters $\lambda$, $\eta$ and $\alpha$ are fixed to $e^8$, 0.8 and 10. Subsampling algorithm will be repeated for 1000 times under each $C$. Mean absolute error, CPU runtime, steps of Newton's Method and their standard deviation are shown in Table \ref{tab:CT}. Considering mean absolute error as the criterion to measure the prediction performance, POSP is always the best of the three types of SSPs. We are not surprised to observe the experimental phenomenon that POSP consumes more CPU runtime than uniform SSP, because the CPU runtime of Algorithm \ref{alg:2step-alg} consists of two parts: time to calculate POSP SSP; time consumed by Newton's Method. Meanwhile BLEV consumes most of CPU runtime on SVD decomposition for full sample matrix. It is also observed that subsampling procedure slightly reduces the iteration steps of Newton's Method. 
\section{Discussions}
We conclude our work in this article from the following open points of view. First, we establish a smooth approximation to the criterion function of Elastic-net. Thereby, a heuristic regression method derived from Elastic-net regression is developed (we call it ‘smooth Elastic-net’). Although the smooth Elastic-net performs well in \S \ref{sec:simulation} and \S \ref{sec:realdata}, its properties such as operational characteristic, grouping effect and solution path need to be examined and studied further.
\par
In addition, POSP in form of \eqref{apposp} is derived from A-optimality which aims at minimizing the asymptotic MSE of the estimator. Nevertheless, the feasibility and performance of SSPs for smooth Elastic-net under other optimal criteria summarized in \citep{morris2010design} like L-optimality, D-optimality, etc., need to be investigated in our future work.
\par
Finally, as a regularization and variable selection method, the Elastic-net is particularly useful when the number of predictors is much bigger than the number of observations, i.e., $p\gg N$. Previously, the conventional approach used to handle sample matrix calculations with equally large $p$ and $N$ was to perform principal component analysis (PCA) on the samples preliminarily \citep[e.g.][]{graf20112d}. The regularization regression (e.g., Elastic-net and Lasso) equipped with subsampling introduces a pioneering approach to address computational complexity problems caused by large-scale sample matrices.
\section*{Supplementary Materials}
The supplementary materials contain the proofs of all the theoretical results in the main text.
\section*{Acknowledgements}
\par

\bibliographystyle{chicago}
\bibliography{myref}






\newpage
\section*{Appendix}

\newpage

\begin{table}[htbp]
    \centering
    \resizebox{\linewidth}{!}{
    \begin{tabular}{r|ccc|ccc|ccc}
          & \multicolumn{3}{c}{mean abs error} \vline & \multicolumn{3}{c}{CPU runtime} \vline& \multicolumn{3}{c}{Newton steps}\\ \hline
         C    &1e3&5e3&1e4&1e3&5e3&1e4&1e3&5e3&1e4 \\\hline
         POSP &6.9949&6.7384&6.7118&0.5844&0.7071&0.8975&14.2420&13.5910&13.4690\\ 
         std  &0.0693&0.0248&0.0170&0.1413&0.1102&0.1165&1.0689&0.9511&0.8921\\\hline
         BLEV &7.0915&6.7587&6.7228&1.2305&0.9389&1.1086&14&13.5070&13.5310 \\
         std &0.1016&0.0363&0.0251&0.9152&0.4407&0.3756&1.0617&0.9396&0.8751\\\hline
         UNIFORM &7.0187&6.7458&6.7158&0.1350&0.2874&0.4853&14.0110&13.5520&13.4950 \\ 
         std &0.0786&0.0307&0.0206&0.0316&0.0560&0.0761&1.0611&0.9457&0.9224\\ \hline
         Full data & \multicolumn{3}{c}{6.6893} \vline & \multicolumn{3}{c}{2.1585} \vline &\multicolumn{3}{c}{15} \\ \hline
    \end{tabular}
    
    }
    
    \caption{2D CT slice experiment result under different $C$, including mean absolute error, CPU runtime, iteration steps of Newton's Method and their standard deviation.}
    \label{tab:CT}
\end{table}
\end{document}